\documentclass[11pt]{article}
\usepackage{amsfonts}
\usepackage{latexsym,amsmath}
\usepackage{amssymb,array}
\usepackage{calc}
\RequirePackage[dvips]{graphicx}
\usepackage{graphics}
\usepackage[colorlinks=true]{hyperref}
\hypersetup{urlcolor=blue, citecolor=red}
\makeatletter \oddsidemargin 0in \evensidemargin 0in \textwidth
16cm\textheight 20cm 
\begin{document}
\newcommand{\la}{\lambda}
\newcommand{\eq}{\Leftrightarrow}
\newcommand{\mf}{\mathbf}
\newcommand{\ri}{\Rightarrow}
\newtheorem{t1}{Theorem}[section]
\newtheorem{d1}{Definition}[section]
\newtheorem{n1}{Notation}[section]
\newtheorem{c1}{Corollary}[section]
\newtheorem{l1}{Lemma}[section]
\newtheorem{r1}{Remark}[section]
\newtheorem{e1}{Counterexample}[section]
\newtheorem{e}{Example}[section]
\newtheorem{re1}{Result}[section]
\newtheorem{p1}{Proposition}[section]
\newtheorem{cn1}{Conclusion}[section]
\renewcommand{\theequation}{\thesection.\arabic{equation}}
\pagenumbering{arabic}
\title {Comparisons of Order Statistics from Some Heterogeneous Discrete Distributions}
\author{Shovan Chowdhury\footnote{Corresponding
author e-mail: shovanc@iimk.ac.in; meetshovan@gmail.com}\;\\Quantitative Methods and Operations Management Area\\Indian Institute of Management, Kozhikode\\Kerala, India.\and
Amarjit Kundu\\Department of
Mathematics\\
Raiganj University\\ West Bengal, India.\and Surja Kanta Mishra\\Department of
Mathematics\\
Raiganj University\\ West Bengal, India.}
\maketitle
\begin{abstract}
In this paper, we compare extreme order statistics through vector majorization arising from heterogeneous Poisson and geometric random variables. These comparisons are carried out with respect to usual stochastic ordering.
 
\end{abstract}
 %{\bf AMS 2010 Subject Classifications}: 62G30, 60E15, 60K10

\section{Introduction}
\setcounter{equation}{0}
\hspace*{0.3in} The notion of stochastic order based on majorization (see Marshall et al.~\cite{Maol}) deals with the diversity of the components of a vector in $\Re^{n}$. In Statistics and related disciplines, the concept is used for the purpose of comparing order statistics (especially minimum and maximum) of two vectors, $\mathbf{a}=(a_1,a_2,\dots,a_n)$ and $\mathbf{b}=(b_1,b_2,\dots,b_n),$ say. Now, $\mathbf{a} $ is said to majorize $\mathbf{b}$ (written as $\mathbf{a}\stackrel{m}{\succeq}\mathbf{b}$) if $\sum_{i=1}^j a_{i:n}\le\sum_{i=1}^j b_{i:n}$, for $j=1,\;2,\;\ldots, n-1$ and $\sum_{i=1}^n a_{i:n}=\sum_{i=1}^n b_{i:n}$, where $a_{i:n}$ denotes the $i$th component in the increasing arrangement of the components of $\mathbf{a}$. The definition of majorization implies that the order statistics of $\mathbf{a}$ are more dispersed than those of $\mathbf{b}$ (although the average is the same for both vectors). The concept of majorization order has been used for the last two decades in many diverse areas including management science, economics, physics, actuarial science, reliability theory and survival analysis. Comparison of smallest and largest order statistics from heterogeneous independent random variables 
following specific continuous distribution function can be found in Dykstra et al.~\cite{dkr11}, Fang and Zhang~\cite{fz}, Zhao and Balakrishnan~(\cite{zb11.2}), %Balakrishnan et al.~\cite{ba1}, 
Torrado and Kochar~\cite{tr11}, Kundu and Chowdhury~(\cite{kun2}), 
%Chowdhury and Kundu~\cite{ch} 
and the references therein. One can find such comparisons under the same set-up for a family of continuous distributions in Khaledi et al.~\cite{kh}, Li et al.~\cite{li2}, Hazra et al.~\cite{ha}, and Kundu and Chowdhury~\cite{kuc}. It is to be noted here that results on stochastic comparison of order statistics are obtained largely for the continuous random variables. \\
\hspace*{0.3in} The field of univariate discrete distributions has been investigated by the researchers extensively in the past seven decades (see Johnson et al.~\cite{jo}). Discrete probability distributions and their applications are common in statistics and related disciplines such as reliability, economics, engineering, insurance, business and management, hydrology, epidemiology, and others. Many well-known discrete distributions, such as binomial, Poisson, negative binomial are studied vastly in the literature.\\
In particular, Poisson distribution is one of the most powerful distributions used in practice for modeling count data. The Poisson distribution is specified by a single parameter which defines both the mean and the variance of the distribution. This characteristic of Poisson distribution makes it an excellent choice for modeling over-dispersed data. For instance, the number of customer service calls per hour is an important metric for assessing the adequacy of customer service staffing in a customer care; the number of losses or claims that occur each year in an insurance firm helps reviewing and evaluating business insurance coverage; number of defects per meter is an essential quality performance indicator in a production line; patient footfall in an out patient department (OPD) of a hospital is an important metric for resource utilization and capacity planning. Poisson distribution can be the best choice for modeling each of the instances as mentioned above.
\\ On the other hand, the geometric distribution is usually employed to model the number of failures in a sequence of independent and identically distributed Bernoulli trials before the first success occurs. Besides its theoretical flavor, geometric distribution has been proved of substantial interest in numerous practical applications, viz. quality control, reliability, psychology, finance, ecology and others. Moreover, the distribution is extensively used in the literature to model discrete failure time. For instance, geometric distribution is appropriate when (i) a piece of equipment operates in cycles and the number of cycles prior to failure is observed, (ii) failures occur only due to incoming shocks and number of rounds fired until failure becomes more crucial than age at failure, (iii) a device is monitored only once per time period (e.g., an hour, a day) and the observation is the number of time periods successfully completed prior to failure of the device.
\\\hspace*{0.3in} The number of papers on discrete order statistics is considerably smaller than on continuous case. For a discussion of order statistics from discrete distributions one can refer to Nagaraja~\cite{na}, Arnold et al.~\cite{ar}, and Dembińska~\cite{dem}. Moreover, a handful number of papers have studied the heterogeneous effect on the distribution properties of discrete order statistics (see Davies and Dembińska~\cite{dav} and Xu and Hu~\cite{xu}). Moreover, comparison of order statistics for heterogenous discrete distributions, based on vector majorization has been rarely studied in the literature. The only paper by Chen et al.~\cite{che} compares the order statistics stochastically from heterogeneous negative binomial random variables. Current work is an attempt to derive some results on ordering properties of extremes for Poisson and geometric distributions. The rest of this paper is organized as follows. In Section 2, we have given the required definitions and some useful lemmas that are used throughout the paper. Results concerning stochastic comparison of the minimum and maximum order statistics are derived in Section~3.\\   
 
\hspace*{0.3 in} Throughout the paper, the word increasing (resp. decreasing) and non decreasing (resp. non increasing) are used interchangeably, and $\mathbb{R^+}$ denotes the set of real numbers $\{x:0<x<\infty\}$. We also write $a\stackrel{sign}{=}b$ to mean that $a$ and $b$ have the same sign. %Further, by $a\stackrel{\rm def}=b$ we mean that $b$ is defined as $a$. 

\section{Preliminaries}
\hspace*{0.3 in} Let $X$ and $Y$ be two discrete random variables  with survival functions $\overline{F}_{X}\left(\cdot\right)$ and $\overline{F}_{Y}\left(\cdot\right)$, distribution functions $F_{X}\left(\cdot\right)$ and $F_{Y}\left(\cdot\right) $, hazard rate function $r_{X}\left(\cdot\right)$ and $r_{Y}\left(\cdot\right)$ and reversed hazard rate functions $\overline{r}_{X}\left(\cdot\right)$ and $\overline{r}_{Y}\left(\cdot\right)$ respectively.
\\\hspace*{0.3 in} In order to compare $X$ and $Y$ different order statistics are used for fair and reasonable comparison. In literature many different kinds of stochastic orders have been developed and studied.
The following well known definitions may be obtained in Shaked and Shanthikumar~\cite{shak1}.
\begin{d1}\label{de1}
Let $X$ and $Y$ be two discrete random variables with respective supports $[0,\infty)$. Then, $X$ is said to be smaller than $Y$ 
\begin{enumerate}
\item[(i)] in usual stochastic (st) order, denoted as $X\leq_{st}Y$, if $\bar F_X(t)\leq \bar F_Y(t)$ for all discrete real number $t\in [0,\infty),$ 
%likelihood ratio (lr) order, denoted as $X\leq_{lr}Y$, if 
%$$\frac{f_Y(t)}{f_X(t)}\;\text{is increasing in} \,t\in(l_X,u_X)\cup(l_Y,u_Y);$$
\item[(ii)] hazard rate (hr) order, denoted as $X\leq_{hr}Y$, if  $r_X(t)\geq r_Y(t)$ for all discrete real number $t\in [0,\infty)$.
 %which can equivalently be written as $r_X(t)\geq r_Y(t)$ for all $t$;
 \item[(iii)] reversed hazard rate (rhr) order, denoted as $X\leq_{rhr}Y$, if  $\overline{r}_X(t)\leq \overline{r}_Y(t)$ for all discrete real number $t\in [0,\infty)$.  
\end{enumerate}
\end{d1}
\hspace*{0.3 in} It is well known that the results on different stochastic orders can be established using majorization order(s). Let $I^n$ denote an $n$-dimensional Euclidean space where $I\subseteq\Re$. Further, let $\mathbf{x}=(x_1,x_2,\dots,x_n)\in I^n$ and $\mathbf{y}=(y_1,y_2,\dots,y_n)\in I^n$ be any two real vectors with $x_{(1)}\le x_{(2)}\le\cdots\le x_{(n)}$ ($x_{[1]}\ge x_{[2]}\ge\cdots\ge x_{[n]}$) being the increasing (decreasing) arrangements of the components of the vector $\mathbf{x}$. The following definitions may be found in Marshall \emph{et al.} \cite{Maol}.
\begin{d1}\label{de12}
%\begin{enumerate}
The vector $\mathbf{x} $ is said to majorize the vector $\mathbf{y} $ (written as $\mathbf{x}\stackrel{m}{\succeq}\mathbf{y}$) if
\begin{equation*}\sum_{i=1}^j x_{[i]}\ge\sum_{i=1}^j y_{[i]},\;j=1,\;2,\;\ldots, n-1,\;\;and\; \;\sum_{i=1}^n x_{[i]}=\sum_{i=1}^n y_{[i]}.\end{equation*}
or equivalently,
\begin{equation*}
\sum_{i=1}^j x_{(i)}\le\sum_{i=1}^j y_{(i)},\;j=1,\;2,\;\ldots, n-1,\;\;and \;\;\sum_{i=1}^n x_{(i)}=\sum_{i=1}^n y_{(i)}.
\end{equation*} 
\end{d1}

\begin{d1}\label{de2}
A function $\psi:I^n\rightarrow\Re$ is said to be Schur-convex (resp. Schur-concave) on $I^n$ if 
\begin{equation*}
\mathbf{x}\stackrel{m}{\succeq}\mathbf{y} \;\text{implies}\;\psi\left(\mathbf{x}\right)\ge (\text{resp. }\le)\;\psi\left(\mathbf{y}\right)\;for\;all\;\mathbf{x},\;\mathbf{y}\in I^n.
\end{equation*}
\end{d1}

%\begin{n1}
%Let us introduce the following notations.
%\begin{enumerate}
%\item[(i)] $\mathcal{D}_{+}=\left\{\left(x_{1},x_2,\ldots,x_{n}\right):x_{1}\geq x_2\geq\ldots\geq x_{n}> 0\right\}$.
%\item[(ii)] $\mathcal{E}_{+}=\left\{\left(x_{1},x_2,\ldots,x_{n}\right):0< x_{1}\leq x_2\leq\ldots\leq x_{n}\right\}$.
%\end{enumerate}
%\end{n1}
Let us introduce the following lemma which will be used in the next section to prove the results.
 \begin{l1}\label{l5}
 $\left(\text{A4 of Marshall \emph{et al.}\cite{Maol}}~,p. 84)\right)$~Let $I\subset\Re$ be an
open interval and let $\varphi : I^n\rightarrow\Re$ be a continuously differentiable function. Then the necessary
and sufficient conditions for $\varphi$ to be Schur-convex (Schur-concave) on $I^n$ are, 
\begin{enumerate}
\item[i)] $\varphi$ is symmetric on $I^n$, and
\item[ii)] for all $\mf{z}\in I^n$ and for all $i\neq j,$ $$\left(z_i-z_j\right)\left[\varphi_{(i)}(\mf{z})-\varphi_{(j)}(\mf{z})\right]\geq (\leq) 0,$$
where $\varphi_{(i)}(\mf{z})=\partial\varphi(\mf{z})/\partial z_i$ denotes the partial derivative of $\varphi$ with respect to its $i$th argument.
\end{enumerate}
\end{l1}
\section{Main Results}
Some results on the comparison of minimum and maximum order statistics are derived in this section. In the first subsection, extreme order statistics from heterogeneous Poisson random variables are compared while in the second subsection, results are derived when each random variable follows geometric distribution.
\subsection{Poisson Distribution}
Let $X$ be a non negative random variable satisfying Poisson distribution with parameter $\mu (> 0)$, having probability mass function (PMF) $$ f_r = P\left(X=r\right)=\frac{e^{-\mu}\mu^{r}}{ r!},~r= 0,1,2,\ldots,$$  and cumulative distribution function (CDF)
$$ F_r= P\left(X\leq r\right) = \sum_{\substack{0\leq k\leq r}}\frac{e^{-\mu}\mu^{k}}{ k!} = \frac{\int\limits_{\mu}^{\infty} e^{-t}t^{r} dt}{r!}.$$
Let $ X_1, X_2, \dots, X_n$ be a sequence of random variable satisfying Poisson distribution with parameters $\mu_1, \mu_2, \dots, \mu_n$ respectively. If $F_{n:n}\left(r\right)$ and $\overline {F}_{1:n}\left(r\right)$ be the distribution and survival function corresponding to the random variable $X_{n:n}$ and $X_{1:n}$ respectively, then clearly

$$F_{n:n}\left(r\right)= \prod_{i=1}^{n} P\left(X\leq r\right)= \prod_{i=1}^{n} \frac{\int\limits_{\mu_i}^{\infty} e^{-t}t^{r} dt}{ r!},~r= 0,1,2,\ldots,$$ and 
$$\overline {F}_{1:n}\left(r\right) =\prod_{i=1}^{n}  P\left(X > r\right)=\prod_{i=1}^{n} \left[1- P\left(X\leq r\right)\right]=\prod_{i=1}^{n} \left[1-\frac{\int\limits_{\mu_i}^{\infty} e^{-t}t^{r} dt}{r!}\right] =\prod_{i=1}^{n} \frac{\int\limits_{0}^{\mu_i} e^{-t}t^{r} dt}{r!}.$$
The following lemma will be used to prove the next theorem.
\begin{l1}
For any positive integer $r$, $$l\left(\mu\right) =\frac{e^{-\mu}\mu^{r}}{\int\limits_{\mu}^{\infty} e^{-t} t^{r}dt }$$ is increasing in $\mu$.
\end{l1}
{\bf Proof:} Assume that $Z$ is a random variable having survival function $H\left(z\right)$, such that $$H\left(z\right)= \frac{\int\limits_{z}^{\infty} e^{-t}t^{r} dt}{r!},~z\in \mathfrak{R_+}.$$
So, if $h\left(z\right)$ is the hazard rate function of the random variable $Z$, then it is clear that $l\left(\mu\right)=h\left(\mu\right).$
Again, from Barlow and Proschan \cite{ba}, the random $Z$ has increasing hazard rate function if the density function $f\left(z\right)=$ $\frac{e^{-z}z^{r}}{r!}$ of the random variable $Z$ is log-concave in $z$.
Now, as for all $z\geq 0$,  $$\log f\left(z\right)= -z+r\log z - \log \left(r!\right),$$ then
$$\frac{\partial}{\partial z} \log f\left(z\right)= -1+\frac{r}{z},$$ and
$$\frac{\partial^2}{\partial z^2} \log f\left(z\right)=-\frac{r}{z^2} \leq 0,$$
giving that $f\left(z\right)$ is log-concave in $z$, then $h\left(z\right)$ is increasing in $z\geq 0$. So $l\left(\mu\right)$ is increasing in $\mu$.\hfill$\diamond$\\
Let, $X_1, X_2, \ldots, X_n$ and $X_1^*, X_2^*, \ldots, X_n^*$ be two sets of random variables having Poisson distribution having parameters $\mu_1,\mu_2,\ldots, \mu_n$ and $\mu_1^*,\mu_2^*,\ldots, \mu_n^*$ respectively. Also let, $\mbox{\boldmath $\mu$}=\left(\mu_1,\mu_2,\ldots, \mu_n\right)$ and $\mbox{\boldmath $\mu^*$}=\left(\mu_1^*,\mu_2^*,\ldots, \mu_n^*\right)$. The following theorem shows that majorized parameter vector of one set of random variables implies greater maximum order statistic than the other, in term of usual stochastic ordering.
 \begin{t1}\label{th1}
 For $i=1,2,\ldots, n$, let $X_i$ and $Y_i$ be two sets of mutually independent random variables each satisfying Poisson distribution with parameters $\mu_i$ and $\mu_i^{*}$ respectively. Then $$\mbox{\boldmath $\mu$}\stackrel{m}\succeq \mbox{\boldmath $\mu^*$} \Rightarrow  X_{n:n}\geq_{st}Y_{n:n}.$$
\end{t1}

{\bf Proof:} If $F_{n:n}(r)$ and $G_{n:n}(r)$ denote the distribution functions of $X_{n:n}$ and $Y_{n:n}$ respectively, then, clearly, for $r= 0,1,2,\ldots$, $$F_{n:n}\left(r\right)= \prod_{k=1}^{n} \frac{\int\limits_{\mu_k}^{\infty} e^{-t}t^{r} dt}{r!}\ \text{and}\ G_{n:n}\left(r\right)= \prod_{k=1}^{n} \frac{\int\limits_{\mu_k^{*}}^{\infty} e^{-t}t^{r} dt}{r!}.$$ Let, 
$$\Psi \left(\mbox{\boldmath $\mu$}\right)= \prod_{k=1}^{n} \frac{\int\limits_{\mu_k}^{\infty} e^{-t}t^{r} dt}{r!}.$$ 
Differentiating $\Psi \left(\mbox{\boldmath $\mu$}\right)$ with respect to $\mu_i$, it can be written that, for $r= 0,1,2,\ldots$,
$$\frac{\partial \Psi\left(\mbox{\boldmath $\mu$}\right)}{\partial \mu_i} = \prod_{k\neq i=1}^{n} \frac{\int\limits_{\mu_k}^{\infty} e^{-t}t^{r} dt}{r!} \left[-e^{-\mu_i}\mu_i^{r}\right].$$
So, for $i \neq j$,
$$\left(\mu_i-\mu_j\right)\left(\frac{\partial \Psi}{\partial \mu_i} - \frac{\partial \Psi}{\partial \mu_j}\right)= \left(\mu_i-\mu_j\right) \Psi \left(\mbox{\boldmath $\mu$}\right)\left[\frac{e^{-\mu_j}\mu_j^{r}}{\int\limits_{\mu_j}^{\infty} e^{-t} t^{r}dt } - \frac{e^{-\mu_i}\mu_i^{r}}{\int\limits_{\mu_i}^{\infty} e^{-t} t^{r}dt }\right].$$
$$=\left(\mu_i-\mu_j\right) \Psi \left(\mbox{\boldmath $\mu$}\right)\left[l\left(\mu_j\right)-l\left(\mu_i\right)\right].$$
Thus, if for $i\leq j$, $\mu_i \geq (\leq) \mu_j$, then $\left(\mu_i-\mu_j\right)\frac{\partial \Psi}{\partial \mu_i} - \frac{\partial \Psi}{\partial \mu_j} \leq 0$, giving $\Psi\left(\mbox{\boldmath $\mu$}\right)$ is $s$-concave in $\mbox{\boldmath $\mu$}$, by Lemma \ref{l5}. So, for all non-negative integer $r$, it can be written that,
$$ \mbox{\boldmath $\mu$}\stackrel{m}\succeq \mbox{\boldmath $\mu^*$}\ \text{implies}\ F_{n:n}\left(r\right)\leq G_{n:n}\left(r\right),$$ proving the result.\hfill$\diamond$\\
Now the question arises: whether the result of the above theorem can be upgraded to reversed hazard rate ordering or not. Next counter-example shows that this is not possible.
\begin{e1} Let, for $n=3$ and for $r=1,2,\ldots$, $\tilde{h}_{n:n}(r)$ and $\tilde{h}^*_{n:n}(r)$ denote the reversed hazard rate functions of $X_{n:n}$ and $Y_{n:n} $ respectively. Now, if $\mbox{\boldmath $\mu$}=(8,0.8,0.1)$ and $\mbox{\boldmath $\mu^*$}=(7,1,0.9)$ are taken, then  although $ \mbox{\boldmath $\mu$}\stackrel{m}\succeq \mbox{\boldmath $\mu^*$}$, it can be observed that for $r=5$, $\tilde{h}_{n:n}(r)-\tilde{h}^*_{n:n}(r)=0.0520158$ and for $r=2$, $\tilde{h}_{n:n}(r)-\tilde{h}^*_{n:n}(r)=-0.0232122$, proving that there exist no reversed hazard rate ordering between $X_{n:n}$ and $Y_{n:n} $.
\end{e1}
%Let $$\Psi_1 \left(\overline{\boldmath \mu}\right)= \prod_{k=1}^{n} \frac{\int\limits_{0}^{\mu_k} \exp^{-t}t^{r} dt}{ \left(r\right)!}$$ 
%$$ \Rightarrow \frac{\partial \Psi_1}{\partial \mu_i} = \prod_{k\neq i=1}^{n} \frac{\int\limits_{0}^{\mu_k} \exp^{-t}t^{r} dt}{ \left(r\right)!} \left[\frac{\exp^{-\mu_i}\mu_i^{r}}{ \left(r\right)!}\right].$$
The next theorem shows that there exits stochastic ordering between $X_{1:n}$ and $Y_{1:n}$, if there exists majorization ordering between between their parameter vectors. 

\begin{t1}\label{th2}
 For $i=1,2,\ldots, n$, let $X_i$ and $Y_i$ be two sets of mutually independent random variables each satisfying Poisson distribution with parameters $\mu_i$ and $\mu_i^{*}$ respectively. Then $$ \mbox{\boldmath $\mu$}\stackrel{m}\succeq \mbox{\boldmath $\mu^*$}\ \text{implies}\ X_{1:n}\leq_{st}Y_{1:n}.$$
\end{t1}
{\bf Proof:} If $\overline{F}_{1:n}(r)$ and $\overline{G}_{1:n}(r)$ denote the survival functions of $X_{1:n}$ and $Y_{1:n}$ respectively, then, clearly, for $r=0,1,2,\ldots$, $$\overline {F}_{1:n}\left(r\right) = \prod_{k=1}^{n} \frac{\int\limits_{0}^{\mu_k} e^{-t}t^{r} dt}{r!}\ \text{and}\ \overline{G}_{1:n}\left(r\right) = \prod_{k=1}^{n} \frac{\int\limits_{0}^{\mu_k^{*}} e^{-t}t^{r} dt}{r!}.$$
Let, for $r=0,1,2,\ldots$ $$\Psi_1 \left(\mbox{\boldmath $\mu$}\right)= \prod_{k=1}^{n} \frac{\int\limits_{0}^{\mu_k} e^{-t}t^{r} dt}{r!}.$$
Differentiating the above expression with respect to $\mu_i$,
$$\frac{\partial \Psi_1\left(\mbox{\boldmath $\mu$}\right)}{\partial \mu_i}=\prod_{k\neq i=1}^{n} \frac{\int\limits_{0}^{\mu_k} e^{-t}t^{r} dt}{r!}\left[e^{-\mu_i}\mu_i^{r}\right].$$
So, for $i \neq j$,
\begin{equation}\label{e0}
\left(\mu_i-\mu_j\right)\left(\frac{\partial \Psi_1\left(\mbox{\boldmath $\mu$}\right)}{\partial \mu_i} - \frac{\partial \Psi_1\left(\mbox{\boldmath $\mu$}\right)}{\partial \mu_j}\right)= \left(\mu_i-\mu_j\right)\Psi_1 \left(\mbox{\boldmath $\mu$}\right)\left[\frac{e^{-\mu_i}\mu_i^{r}}{\int\limits_{0}^{\mu_i} e^{-t} t^{r}dt } - \frac{e^{-\mu_j}\mu_j^{r}}{\int\limits_{0}^{\mu_j} e^{-t} t^{r}dt }\right].
\end{equation}
Now, it can be noted that, for $r=0,1,2,\ldots$, $$\frac{e^{-\mu}\mu^{r}}{\int\limits_{0}^{\mu} e^{-t} t^{r}dt }= \frac{1}{\int\limits_{0}^{\mu} e^{\mu-t} \left(\frac{t}{\mu}\right)^{r}dt }.$$
So, taking $ 1-\frac{t}{\mu}=u$, the above expression can be rewritten as $$\frac{e^{-\mu}\mu^{r}}{\int\limits_{0}^{\mu} e^{-t} t^{r}dt }= \frac{1}{-\int\limits_{1}^{0} e^{\mu u} \left( 1-u\right)^{r} \mu du  }=  \frac{1}{\mu \int\limits_{0}^{1} e^{\mu u} \left( 1-u\right)^{r} du  }.$$
Now, as for $r=0,1,2,\ldots$, 
$$\frac{\mathrm d}{\mathrm d \mu} \left( \int\limits_{0}^{1} e^{\mu u} \left( 1-u\right)^{r} du \right)= \int\limits_{0}^{1} u e^{\mu u} \left( 1-u\right)^{r} du  \geq 0,$$ then $\int\limits_{0}^{1} e^{\mu u} \left( 1-u\right)^{r} du$ is increasing in $\mu$. So, by noticing the fact that $ \int\limits_{0}^{1} e^{\mu u} \left( 1-u\right)^{r} du \geq 0$, it can be written that $\mu \int\limits_{0}^{1} e^{\mu u} \left( 1-u\right)^{r} du$ is also increasing in $\mu$. So, for $i\leq j$, $\mu_i \geq (\leq) \mu_j$ implies that for $r=0,1,2,\ldots$,
$$\frac{1}{\mu_i \int\limits_{0}^{1} e^{\mu_i u} \left( 1-u\right)^{r}  du} \leq (\geq)\frac{1}{\mu_j \int\limits_{0}^{1} e^{\mu_j u} \left( 1-u\right)^{r}  du},$$
 which gives
$$\frac{e^{-\mu_i}\mu_i^{r}}{\int\limits_{0}^{\mu_i} e^{-t} t^{r}dt} \leq (\geq)\frac{e^{-\mu_j}\mu_j^{r}}{\int\limits_{0}^{\mu_j} e^{-t} t^{r}dt}.$$
So, from (\ref{e0}), it can be written that $\left(\mu_i-\mu_j\right)\left(\frac{\partial \Psi_1}{\partial \mu_i} - \frac{\partial \Psi_1}{\partial \mu_j}\right) \leq 0$, giving $\Psi_1\left(\mbox{\boldmath $\mu$}\right)$ is $s$-concave in $\mbox{\boldmath $\mu$}$, by Lemma \ref{l5}. So, for all $r=0,1,2,\ldots$, $$ \mbox{\boldmath $\mu$}\stackrel{m}\succeq \mbox{\boldmath $\mu^*$}\ \text{implies}\ \overline{F}_{1:n}(r)\leq \overline{G}_{1:n}(r),$$ proving the result.\hfill$\diamond$\\
Next one counter-example is given to show that the above theorem cannot be improved further to hazard rate ordering. 
\begin{e1} For $n=3$ and for $r=1,2,\ldots$, let $h_{n:n}(r)$ and $h^*_{n:n}(r)$ be the hazard rate functions of $X_{1:n}$ and $Y_{1:n} $ respectively. Now, if $\mbox{\boldmath $\mu$}=(28,0.8,0.1)$ and $\mbox{\boldmath $\mu^*$}=(27,1,0.9)$ are taken, then although $ \mbox{\boldmath $\mu$}\stackrel{m}\succeq \mbox{\boldmath $\mu^*$}$, it can be observed that for $r=16$, $h_{1:n}(r)-h^*_{1:n}(r)=-0.00024431$ and for $r=6$, $h_{1:n}(r)-h^*_{n:n}(r)=0.0124328$, proving that there exist no hazard rate ordering between $X_{1:n}$ and $Y_{1:n} $.
\end{e1}

\subsection{Geometric Distribution}
Let $ X_1, X_2, \ldots, X_n$ be a sequence of random variables satisfying geometric distribution with parameter $q_1, q_2, \ldots, q_n$. Now, if  $\overline {F}_{1:n}\left(u\right)$ and $F_{n:n}\left(u\right)$ be the survival function and distribution functions of $X_{1:n}$ and $X_{n:n}$ respectively, then, for all $ u= 0,1,2,\ldots$,
$$ \overline {F}_{1:n}\left(u\right) = \prod_{i=1}^{n} q_i^{u+1}\ \text{and}\ F_{n:n}\left(u\right)=\prod_{k=1}^{n} \left(1-q_k^{u+1}\right).$$ 
Now, if $r_{1:n}\left(u\right)$ denotes hazard rate function of the random variable $X_{1:n}$, then, for all $ u= 0,1,2,\ldots$, it can be written that,
$$ r_{1:n}\left(u\right)= \frac{\overline {F}_{1:n}\left(u\right)-\overline {F}_{1:n}\left(u+1\right)}{\overline {F}_{1:n}\left(u\right)}= 1- \frac{\prod_{i=1}^{n} q_i^{u+1}}{\prod_{i=1}^{n} q_i^{u}}= 1-\prod_{i=1}^{n} q_i.$$
Let, $X_1, X_2, \ldots, X_n$ and $X_1^*, X_2^*, \ldots, X_n^*$ be two sets of random variables following geometric distribution having parameters $q_1,q_2,\ldots, q_n$ and $q_1^*,q_2^*,\ldots, q_n^*$ respectively. Also let, $\mbox{\boldmath $q$}=\left(q_1,q_2,\ldots, q_n\right)$ and $\mbox{\boldmath $q^*$}=\left(q_1^*,q_2^*,\ldots, q_n^*\right)$. The following theorem shows that, under certain condition $X_{1:n}$ is greater than $Y_{1:n}$ in hazard rate ordering. The proof is straight forward and hence omitted.
\begin{t1}\label{th3}
For $i=1,2,\ldots, n$, let $X_i$ and $Y_i$ be two sets of mutually independent random variables each satisfying geometric distribution with parameters $q_i$ and $q_i^{*}$ respectively. Then $$ \prod_{i=1}^{n} q_i\geq \prod_{i=1}^{n} q_i^{*}\ \text{implies}\  X_{1:n}\geq_{hr}Y_{1:n}.$$
\end{t1}
The following theorem shows that under certain condition $X_{n:n}$ will be greater than $Y_{n:n}$ in usual stochastic ordering.
\begin{t1}\label{th4}
For $i=1,2,\ldots, n$, let $X_i$ and $Y_i$ be two sets of mutually independent random variables each satisfying geometric distribution with parameters $q_i$ and $q_i^{*}$ respectively. Then $$\mbox{\boldmath $q$}\stackrel{m}\succeq \mbox{\boldmath $q^*$}\ \text{implies}\ X_{n:n}\geq_{st}Y_{n:n}$$
\end{t1}
{\bf Proof:} If $F_{n:n}\left(u\right)$ and $G_{n:n}\left(u\right)$ be the distribution functions of $X_{n:n}$ and $Y_{n:n}$ respectively, then clearly, for $u=0,1,2,\ldots$,
$$F_{n:n}\left(u\right)=\prod_{k=1}^{n} \left(1-q_k^{u+1}\right)\ \text{and}\ G_{n:n}\left(u\right)=\prod_{k=1}^{n} \left(1-q_k^{* u+1}\right).$$  
Let, $$\Psi_2 \left(\mbox{\boldmath $q$}\right)= \prod_{k=1}^{n} \left(1-q_k^{u+1}\right).$$ 
Differentiating the above equation with respect to $q_i$, it can be written that
$$\frac{\partial \Psi_2\left(\mbox{\boldmath $q$}\right)}{\partial q_i} = \prod_{k\neq i=1}^{n} \left(1-q_k^{u+1}\right)\left[-\left(u+1\right)q_i^{u}\right].$$
So, for $i \leq j$ and for all $u=0,1,2,\ldots$,
\begin{eqnarray*}
\left(q_i-q_j\right)\left(\frac{\partial \Psi_2}{\partial q_i} - \frac{\partial \Psi_2}{\partial q_j}\right)&=& \left(q_i-q_j\right) \left[\prod_{k\neq i,j=1}^{n} \left(1-q_k^{u+1}\right)\right]\left[\left(u+1\right)q_j^{u}\left(1-q_i^{u+1}\right) - \left(u+1\right)q_i^{u}\left(1-q_j^{u+1}\right)\right]\\
 &\stackrel{sign}{=}& \left(q_i-q_j\right)\left(q_j^{u}-q_i^{u}+ q_j^{u}q_i^{u}(q_j-q_i)\right).
\end{eqnarray*}
So, for $i\leq j$, $q_i \geq (\leq) q_j$ gives $\left(q_i-q_j\right)\left( \frac{\partial \Psi_2}{\partial q_i} - \frac{\partial \Psi_2}{\partial q_j}\right) \leq 0$, giving $\Psi_2\left(\mbox{\boldmath $q$}\right)$ is $s$-concave in $\mbox{\boldmath $q$}$, by Lemma \ref{l5}. Thus
$$ \mbox{\boldmath $q$}\stackrel{m}\succeq \mbox{\boldmath $q^*$}\ \text{implies}\ F_{n:n}(u)\leq G_{n:n}(u),$$
proving the result.\hfill$\diamond$\\

Although there exists stochastic ordering between $X_{n:n}$ and $Y_{n:n}$ when $ \mbox{\boldmath $q$}\stackrel{m}\succeq \mbox{\boldmath $q^*$}$, the following counterexample shows that there exists no reversed hazard rate ordering between them. 
\begin{e1}
Let, for $n=3$ and for $r=1,2,\ldots$, $\tilde{h}_{n:n}(r)$ and $\tilde{h}^*_{n:n}(r)$ denote the reversed hazard rate functions of $X_{n:n}$ and $Y_{n:n} $ respectively. Now, if $\mbox{\boldmath $q$}=(0.99,0.96,0.57)$ and $\mbox{\boldmath $q^*$}=(0.9,0.78,0.57)$ are taken, then although $ \mbox{\boldmath $q$}\stackrel{m}\succeq \mbox{\boldmath $q^*$}$, it can be observed that for $u=1$, $\tilde{h}_{n:n}(u)-\tilde{h}^*_{n:n}(u)=-0.0010584$ and for $u=4$, $\tilde{h}_{n:n}(u)-\tilde{h}^*_{n:n}(u)=0.00628996$, proving that there exist no reversed hazard rate ordering between $X_{n:n}$ and $Y_{n:n} $.
\end{e1}

\end{document}